# EFFICIENT AND ACCURATE NUMERICAL SCHEMES FOR A HYDRODYNAMICALLY COUPLED PHASE FIELD DIBLOCK COPOLYMER MODEL

QING CHENG [†], XIAOFENG YANG [‡], JIE SHEN [§*]

ABSTRACT. In this paper, we consider numerical approximations of a hydrodynamically coupled phase field diblock copolymer model, in which the free energy contains a kinetic potential, a gradient entropy, a Ginzburg-Landau double well potential, and a long range nonlocal type potential. We develop a set of second order time marching schemes for this system using the "Invariant Energy Quadratization" approach for the double well potential, the projection method for the Navier-Stokes equation, and a subtle implicit-explicit treatment for the stress and convective term. The resulting schemes are linear and lead to symmetric positive definite systems at each time step, thus they can be efficiently solved. We further prove that these schemes are unconditionally energy stable. Various numerical experiments are performed to validate the accuracy and energy stability of the proposed schemes.

## 1. INTRODUCTION

Block copolymer is a linear-chain molecule composed of two or more subchains linked together to create a polymer chain. When the subchain is made of two (or three) distinct monomer blocks, it is called diblock (or triblock) copolymer. Due to the incompatibility between the blocks, block copolymers undergo a micro-phase separation to form a periodic morphology in nanoscale, which provides an efficient technique to produce nano-structured materials and nano-devices (cf. [12, 19, 20, 29, 31, 44]).

Modeling and numerical simulation are effective means to investigate the phase separation behaviors of block copolymers. In this paper, we consider the phase field based model for diblock copolymer (PF-BCP) model (cf. [4, 7, 8, 16, 24, 30, 44, 46]) known as the dynamic mean field theory, where an order parameter is used to denote the difference between the local volume fractions of two monomers. The evolution of the PF-BCP system is derived from the energetic variation of the action function of the total free energy in the $H^{-1}$ Sobolev space, i.e., the Cahn-Hilliard (CH) type equation. The total free energy for the system is the nonlocal Ohta-Kawasaki functional, that is the standard Cahn-Hilliard free energy supplemented with a nonlocal term, reflecting the first order effects of the connectivity of the monomer chains (cf. [7, 30]).

From a numerical point of view, for the model with stiffness, it is advantageous to develop schemes with unconditionally energy stability in both semi-discrete and fully discrete cases. The main difficulty in designing such schemes is how to discretize those nonlinear terms that involve the stiffness issue originated from the thin interface parameter. In fact, the simple fully-implicit or explicit type discretizations will induce very severe time step constraint depending on the interfacial width (cf. [3, 13, 35]), so they are not efficient in practice. There have been some efforts on developing energy stable schemes for the PF-BCP model recently (cf. [1, 23, 44]). These schemes are essentially based on either the nonlinear convex splitting approach (cf. [1, 11]) or the linear stabilization approach (cf. [6, 26, 28, 33–40, 45, 47, 49, 50, 54–56]), and only first order accurate in time. Undoubtedly, higher order time marching schemes are preferable to lower order schemes if the adopted time step is expected to be as large as possible under certain accuracy requests. This fact motivates us to develop more accurate schemes, e.g., the second order time marching schemes while preserving the unconditional energy stability in this paper.

[†]School of Mathematical Sciences and Fujian Provincial Key Laboratory on Mathematical Modeling and High Performance Scientific Computing, Xiamen University, Xiamen, Fujian, P. R. China, 361005; Email: chengqing@stu.xmu.edu.cn

[‡]Department of Mathematics, University of South Carolina, Columbia, SC, 29208; Email: xfyang@math.sc.edu.

[§*]Corresponding author, School of Mathematical Sciences and Fujian Provincial Key Laboratory on Mathematical Modeling & High Performance Scientific Computing, Xiamen University, Xiamen, Fujian 361005, P. R. China; Department of Mathematics, Purdue University, West Lafayette, IN, USA, 47907; Email:shen7@purdue.edu.

Instead of using traditional discretization approaches like simple implicit, stabilized explicit, convex splitting, or other various tricky Taylor expansions to discretize the nonlinear potentials, we adopt the so-called *Invariant Energy Quadratization* (IEQ) method, which is a novel approach developed recently and successfully applied for a number of gradient flow type models (cf. [21, 48, 51, 52]). The essential idea of the IEQ method is to transform the free energy into a quadratic form (since the nonlinear potential is usually bounded from below) of a set of new variables via a change of variables. The new, equivalent system still retains the similar energy dissipation law in terms of the new variables. For the time-continuous case, the energy law of the new reformulated system is equivalent to the energy law of the original system. One great advantage of such reformulation is that all nonlinear terms can be treated semi-explicitly in a time marching scheme, leading to a linear system at each time step. Moreover, the resulted linear system is symmetric positive definite, thus it can be efficiently solved by a Krylov subspace method such as CG.

Based on this new approach, we develop in this paper a set of efficient schemes which are *accurate* (second order in time), *easy-to-implement* (linear), and *unconditionally energy stable* (with a discrete energy dissipation law) to solve the PF-BCP model and the PF-BCP model coupled with Navier-Stokes equations (PF-BCP-NS). We use the projection method to solve Navier-Stokes equations, and a subtle implicit-explicit treatment to treat the convective and stress terms. We rigorously prove that the unconditionally energy stabilities hold for all proposed schemes. We demonstrate the stability and the accuracy of the proposed schemes through various classical benchmark simulations. To the best of the authors' knowledge, the proposed schemes are the first second order accurate schemes for the PF-BCP model with unconditional energy stabilities.

The rest of the paper is organized as follows. In Section 2, we describe the BCP model with and without hydrodynamics field. In Section 3, we develop the numerical schemes and prove the well-posedness of the linear system, as well as the unconditional energy stabilities. In Section 4, we we present various numerical simulations to validate our numerical schemes. Finally, some concluding remarks are presented in Section 5.

## 2. Model Equations

Let us first introduce some notations. We denote by $(f,g) = (\int_\Omega f(\boldsymbol{x})g(\boldsymbol{x})d\boldsymbol{x})^{\frac{1}{2}}$ the $L^2$ inner product between functions $f(\boldsymbol{x})$ and $g(\boldsymbol{x})$, by $\|f\| = (f,f)$ the $L^2$ norm of function $f(\boldsymbol{x})$. For any $\phi \in L^2(\Omega)$, we denote $\overline{\phi} = \frac{1}{|\Omega|} \int_\Omega \phi(\boldsymbol{x})d\boldsymbol{x}$. Let us define

$$(2.1) \qquad L_0^2(\Omega) = \{\phi \in L^2(\Omega) : \int_\Omega \phi d\boldsymbol{x} = 0\},$$

and the inverse Laplace operator $(-\Delta)^{-1}$: $u \in L_0^2(\Omega) \mapsto v := (-\Delta)^{-1}u$ by

$$(2.2) \quad \begin{cases} -\Delta v = u, \quad \int_\Omega v d\boldsymbol{x} = 0, \\ \text{with the boundary conditions either } (i) \ v \text{ is periodic, or } (ii) \ \partial_{\mathbf{n}} v|_{\partial\Omega} = 0. \end{cases}$$

We now give a brief introduction to the PF-BCP model. The state of the system is described by the local volume fraction of two monomers, $\phi(\boldsymbol{x}, t)$, at all points $\boldsymbol{x} \in \Omega^d, d = 2, 3$ and at time $t$. The total phenomenological free energy is as follows [1, 4, 7, 8, 16, 23, 24, 30, 44, 44, 46]:

$$(2.3) \quad E(\phi) = \int_\Omega \Big(\frac{\epsilon^2}{2}|\nabla\phi|^2 + F(\phi)\Big)d\boldsymbol{x} + \frac{\alpha}{2}\int_\Omega\int_\Omega G(\boldsymbol{x}-\boldsymbol{y})(\phi(\boldsymbol{x})-\overline{\phi})(\phi(\boldsymbol{y})-\overline{\phi})d\boldsymbol{y}d\boldsymbol{x},$$

where the $F(\phi) = \frac{1}{4}(\phi^2-1)^2$ is the Ginzburg-Landau double well potential, $\epsilon$ is the gradient energy coefficient, $\alpha$ is a positive phenomenological parameter, $G$ denotes the Green's function such that $\Delta G(\boldsymbol{x}-\boldsymbol{y}) = -\delta(\boldsymbol{x}-\boldsymbol{y})$ with periodic boundary condition, and $\delta$ is a Dirac delta function. It is clear that the energy functional $E(\phi)$ is the commonly used Cahn-Hilliard free energy when $\alpha = 0$. With the nonlocal term when $\alpha \neq 0$, $E(\phi)$ is referred to as the Ohta-Kawasaki functional, that was first proposed in [30].

The phenomenological mesoscopic dynamic equation is the Cahn-Hilliard type, i.e., the gradient flow system in $H^{-1}(\Omega)$:

$$(2.4) \quad \phi_t = M(\Delta\mu - \alpha(\phi - \overline{\phi})),$$
$$(2.5) \quad \mu = -\epsilon^2\Delta\phi + f(\phi),$$

where $M$ is the mobility constant, $f(\phi) = F'(\phi) = \phi(\phi^2-1)$. The boundary conditions can be either one of the following two type:

$$(2.6) \quad (i) \ \phi, \mu \text{ are periodic; or } (ii) \ \partial_{\mathbf{n}}\phi|_{\partial\Omega} = \nabla\mu \cdot \mathbf{n}|_{\partial\Omega} = 0,$$

where $\mathbf{n}$ is the unit outward normal on the boundary $\partial\Omega$.

In order to derive an energy law, we shall first reformulate the system (2.4)-(2.5) into an equivalent system for which an energy law can be easily derived.

Since $\int_\Omega (\phi - \overline{\phi})d\boldsymbol{x} = 0$, thus we can set

$$(2.7) \quad \psi = (-\Delta)^{-1}(\phi - \overline{\phi}).$$

Let $w = \mu + \alpha\psi$, then the system (2.4)-(2.5) can be rewritten as

$$(2.8) \quad \phi_t = M\Delta w,$$
$$(2.9) \quad w = -\epsilon^2\Delta\phi + f(\phi) + \alpha\psi,$$

with either boundary condition from

$$(2.10) \quad (i) \ \phi, w \text{ are periodic; or } (ii) \ \partial_{\mathbf{n}}\phi|_{\partial\Omega} = \nabla w \cdot \mathbf{n}|_{\partial\Omega} = 0.$$

**Lemma 2.1.** *The system* (2.8)-(2.9) *with* (2.10) *admits the following energy law:*

$$(2.11) \quad \frac{d}{dt}E(\phi) = -M\|\nabla w\|^2,$$

*where*

$$(2.12) \quad E(\phi) = \int_\Omega \Big(\frac{\epsilon^2}{2}|\nabla\phi|^2 + F(\phi) + \frac{\alpha}{2}|\nabla\psi|^2\Big)d\boldsymbol{x},$$

*and $\psi$ is given by* (2.7).

*Proof.* By taking the $L^2$ inner product of (2.8) with $-w$, and of (2.9) with $\phi_t$ and performing integration by parts, we obtain

$$-(\phi_t, w) = M\|\nabla w\|^2, \tag{2.13}$$

$$(w, \phi_t) = \frac{d}{dt}\int_\Omega \Big(\frac{\epsilon^2}{2}|\nabla\phi|^2 + F(\phi)\Big)d\bm{x} + \alpha(\psi, \phi_t). \tag{2.14}$$

We derive from the definition of $\psi$ in (2.7) that

$$-\Delta\psi_t = \phi_t - \frac{1}{|\Omega|}\int_\Omega \phi_t d\bm{x}. \tag{2.15}$$

By taking the $L^2$ inner product of (2.15) with $\alpha\psi$ and notice $(\psi, 1) = \int_\Omega \psi d\bm{x} = 0$, we obtain

$$\frac{d}{dt}\int_\Omega \frac{\alpha}{2}|\nabla\psi|^2 d\bm{x} = \alpha(\phi_t, \psi). \tag{2.16}$$

By combining (2.13), (2.14) and (2.16), we obtain (2.11). $\square$

It is clear that the energy functional (2.12) is equivalent to the energy functional defined in (2.3), and that the potential $w$ is the variational derivative of $E(\phi)$, i.e., $w = \frac{\delta E(\phi)}{\delta \phi}$. By taking the $L^2$ inner product with 1 for (2.8), we obtain the mass conservation property as

$$\frac{d}{dt}\int_\Omega \phi d\bm{x} = 0. \tag{2.17}$$

When coupling with the flow field, the total free energy becomes

$$E(\phi, \bm{u}) = \int_\Omega \Big(\frac{1}{2}|\bm{u}|^2 + \lambda\big(\frac{\epsilon^2}{2}|\nabla\phi|^2 + F(\phi) + \frac{\alpha}{2}|\nabla\psi|^2\big)\Big)d\bm{x}. \tag{2.18}$$

where $\bm{u}$ is the fluid velocity, $\lambda$ is the magnitude of the free energy. Assuming the fluid is incompressible, and following a generalized Fick's law that the mass flux be proportional to the gradient of the chemical potential [5, 14, 25, 27], we can derive the following PF-BCP-NS model:

$$\phi_t + \nabla\cdot(\bm{u}\phi) = M\Delta w, \tag{2.19}$$

$$w = \lambda(-\epsilon^2\Delta\phi + f(\phi) + \alpha\psi), \tag{2.20}$$

$$\bm{u}_t + (\bm{u}\cdot\nabla)\bm{u} + \nabla p - \nu\Delta\bm{u} + \phi\nabla w = 0, \tag{2.21}$$

$$\nabla\cdot\bm{u} = 0, \tag{2.22}$$

where $p$ is the pressure, $\nu$ is the viscosity and $\psi = (-\Delta)^{-1}(\phi - \overline{\phi})$.

The boundary conditions can be either one of the following two types:

$$(i)\ \bm{u}, \phi, w\ \text{are periodic; or}\ (ii)\ \bm{u}|_{\partial\Omega} = \partial_{\bm{n}}\phi|_{\partial\Omega} = \nabla w\cdot\bm{n}|_{\partial\Omega} = 0. \tag{2.23}$$

**Lemma 2.2.** *The system* (2.19)*-*(2.22) *with* (2.23) *admits the following energy law:*

$$\frac{d}{dt}E(\phi, \bm{u}) = -M\|\nabla w\|^2 - \nu\|\nabla\bm{u}\|^2 \leq 0. \tag{2.24}$$

*Proof.* By taking the $L^2$ inner product of (2.19) with $w$, of (2.20) with $-\phi_t$, of (2.21) with $\bm{u}$, and adding the results together, we can obtain (2.24).

$\square$

## 3. Numerical Schemes

We now develop a set of second order semi-discrete numerical schemes to solve the system (2.8)-(2.9) without flow, and the system (2.19)-(2.22) with flow. While we consider only time discretizations here, the results can carry over to any consistent finite-dimensional Galerkin approximations (finite elements or spectral) since the proof is based on variational formulations with all test functions in the same space as the trial function.

3.1. **PF-BCP model.** The main difficult issue here is how to discretize the nonlinear term $f(\phi)$. We shall handle this term using the IEQ approach [21,48,51,52]. The main idea of the IEQ method is to transform the free energy and the PDE system into equivalent forms in terms of new variables via the change of variables. Thus the nonlinear terms can be treated semi-explicitly.

We define a new variable

$$U = \phi^2 - 1, \tag{3.1}$$

so the total energy (2.12) becomes

$$E(\phi, U) = \int_\Omega \Big(\frac{\epsilon^2}{2}|\nabla\phi|^2 + \frac{1}{4}U^2 + \frac{\alpha}{2}|\nabla\psi|^2\Big)d\boldsymbol{x}, \tag{3.2}$$

and the system (2.8)-(2.9) becomes

$$\phi_t = M\Delta w, \tag{3.3}$$

$$w = -\epsilon^2\Delta\phi + \phi U + \alpha\psi, \tag{3.4}$$

$$U_t = 2\phi\phi_t, \tag{3.5}$$

with initial conditions

$$\phi|_{t=0} = \phi^0, \quad U|_{t=0} = (\phi^0)^2 - 1. \tag{3.6}$$

The boundary conditions are still (2.10) since the equation (3.5) for the new variable $U$ does not involve spatial derivative.

The transformed PDE system (3.3)–(3.5) also admits an energy dissipative law. Indeed, by taking the $L^2$ inner product of (3.3) with $-w$, of (3.4) with $\phi_t$, of (3.5) with $-\frac{1}{2}U$, performing integration by parts and summing up all equalities, we can obtain (3.3)–(3.5) as

$$\frac{d}{dt}E(\phi, U) = -M\|\nabla w\|^2 \leq 0. \tag{3.7}$$

**Remark 3.1.** *The new transformed system (3.3)–(3.5) is equivalent to the original system (2.8)-(2.9) since (3.1) can be easily obtained by integrating (3.5) with respect to the time. Therefore, the energy law (3.7) for the transformed system is exactly the same as the energy law (2.11) for the original system.*

We first construct a second order Crank-Nicolson scheme for (3.3)-(3.5).

Let $\delta t > 0$ be a time step size and set $t^n = n\delta t$ for $0 \leq n \leq N = [T/\delta t]$, let $S^n$ denotes the numerical approximation to $S(\cdot, t)|_{t=t^n}$, and $S^{n+\frac{1}{2}} = \frac{S^{n+1}+S^n}{2}$ for any function $S$.

**Scheme 1.** *Assuming that $\phi^{n-1}$, $\phi^n$ and $U^{n-1}$, $U^n$ are known, we solve $\phi^{n+1}$, $U^{n+1}$ as follows:*

$$\frac{\phi^{n+1} - \phi^n}{\delta t} = M\Delta w^{n+\frac{1}{2}}, \tag{3.8}$$

$$w^{n+\frac{1}{2}} = -\epsilon^2\Delta\phi^{n+\frac{1}{2}} + \phi^{\star,n+\frac{1}{2}}U^{n+\frac{1}{2}} - \alpha\psi^{n+\frac{1}{2}}, \tag{3.9}$$

$$U^{n+1} - U^n = 2\phi^{\star,n+\frac{1}{2}}(\phi^{n+1} - \phi^n). \tag{3.10}$$

where $\phi^{\star,n+\frac{1}{2}} = \frac{3}{2}\phi^n - \frac{1}{2}\phi^{n-1}$, and

$$\psi^{n+\frac{1}{2}} = (-\Delta)^{-1}(\phi^{n+\frac{1}{2}} - \overline{\phi}^{n+\frac{1}{2}}). \tag{3.11}$$

*The boundary conditions are either*

(3.12) $\qquad (i)\ \phi^{n+1}, w^{n+\frac{1}{2}}$ *are periodic; or* $(ii)\ \partial_{\mathbf{n}}\phi^{n+1}|_{\partial\Omega} = \nabla w^{n+\frac{1}{2}} \cdot \mathbf{n}|_{\partial\Omega} = 0.$

*To start off, we can compute $\phi^1$ and $U^1$ by using a first-order version of the above scheme.*

**Remark 3.2.** *By taking the $L^2$ inner product of (3.8) with 1, we obtain*

(3.13) $$\int_\Omega \phi^{n+1} d\boldsymbol{x} = \int_\Omega \phi^n d\boldsymbol{x} = \cdots = \int_\Omega \phi^0 d\boldsymbol{x}.$$

**Remark 3.3.** *We can eliminate $U^{n+1}$ and $w^{n+1}$ from (3.8)-(3.10) to obtain*

(3.14) $$\frac{\phi^{n+1} - \phi^n}{\delta t} = M\left(-\epsilon^2 \Delta^2 \phi^{n+\frac{1}{2}} + \Delta((\phi^{\star,n+\frac{1}{2}})^2 \phi^{n+\frac{1}{2}}) - \alpha(\phi^{n+\frac{1}{2}} - \overline{\phi}^{n+\frac{1}{2}})\right) + g^n,$$

*where $g^n$ contains some explicit terms. We now describe how to solve (3.14) in practice.*

*Applying $(-\Delta)^{-1}$ to the above, we obtain*

(3.15) $$(\frac{2}{M\delta t} + \alpha)(-\Delta)^{-1}\phi^{n+1} - \epsilon^2 \Delta \phi^{n+1} + (\phi^{\star,n+\frac{1}{2}})^2 \phi^{n+1} = h^n,$$

*where $h^n$ includes all explicit terms. Since $(-\Delta)^{-1}$ is an self-adjoint positive definite operator, we derive immediately that the equation (3.15) admits a unique solution.*

*Since $(-\Delta)^{-1}$ is a non-local operator, it is not efficient to implement the scheme in physical space. However, in the phase space formed by the eigenfunctions of Laplace operator, $(-\Delta)^{-1}$ is a local operator, and can be implemented efficiently. More precisely, let $X_N$ be a suitable approximation space, we consider the following Galerkin approximation for (3.15):*

*Find $\phi_N^{n+1} \in X_N$ such that, $\forall v_N \in X_N$, we have*

(3.16) $$(\frac{2}{M\delta t} + \alpha)((-\Delta)^{-1}\phi_N^{n+1}, v_N) - \epsilon^2(\Delta \phi^{n+1}, v_N) + ((\phi^{\star,n+\frac{1}{2}})^2 \phi^{n+1}, v_N) = (h^n, v_N).$$

*It is clear that the above linear system is symmetric positive definite.*

- *For periodic boundary conditions, we choose $X_N$ to be the trigonometric functions, of degree less than or equal to $N$, which are eigenfunctions of Laplace operator are Fourier series so that $(-\Delta)^{-1}\phi_N$ and $\Delta\phi_N$ can be evaluated easily for $\phi_N \in X_N$.*
- *For non-periodic problems, we can use the Fourier-like basis functions [41], which are eigenfunctions of discrete Laplacian operator, so that $(-\Delta)^{-1}\phi_N$ and $\Delta\phi_N$ can be evaluated as in period case.*
- *However, the term $((\phi^{\star,n+\frac{1}{2}})^2 \phi^{n+1}, v_N)$ will lead to a full matrix so a direct solver is not practical. We shall solve (3.16) by using a preconditioned conjugate gradient (PCG) method with an optimal preconditioner contructed by an approximate problem of (3.16) where $\phi^*$ is replaced by a constant $\overline{\phi}^* \approx \int_\Omega \phi^*$. Since we use the Fourier basis (in the period case) or the Fourier-like basis in the non-periodic case, the linear system for the preconditioner is diagonal. Therefore, the system (3.16) can be solved very efficiently.*

As for the energy stability, we have the following theorem.

**Theorem 3.1.** *The scheme (3.8)-(3.10) admits a unique solution, and is unconditionally energy stable, i.e., satisfies the following discrete energy dissipation law:*

(3.17) $$\frac{1}{\delta t}(E_{cn2}(\phi^{n+1}, U^{n+1}) - E_{cn2}(\phi^n, U^n)) = -M\|\nabla w^{n+\frac{1}{2}}\|^2,$$

*where $E_{cn2}(\phi, U) = \frac{\epsilon^2}{2}\|\nabla\phi\|^2 + \frac{1}{4}\|U\|^2 + \frac{\alpha}{2}\|\nabla\psi\|^2$, and $\psi$ is determined by $\phi$ from (3.11).*

*Proof.* Applying $(-\Delta)^{-1}$ to (3.14), by taking the $L^2$ inner product of (3.8) with $-\delta t w^{n+\frac{1}{2}}$, we obtain

$$-(\phi^{n+1} - \phi^n, w^{n+\frac{1}{2}}) = \delta t M \|\nabla w^{n+\frac{1}{2}}\|^2. \tag{3.18}$$

By taking the $L^2$ inner product of (3.9) with $\phi^{n+1} - \phi^n$ and applying the following inequality,

$$2(a-b, a) = |a|^2 - |b|^2 + |a+b|^2, \tag{3.19}$$

we obtain

$$(w^{n+\frac{1}{2}}, \phi^{n+1} - \phi^n) = \frac{\epsilon^2}{2}(\|\nabla \phi^{n+1}\|^2 - \|\nabla \phi^n\|^2) + (\phi^{\star,n+\frac{1}{2}} U^{n+\frac{1}{2}}, \phi^{n+1} - \phi^n) \\ + \alpha(\psi^{n+\frac{1}{2}}, \phi^{n+1} - \phi^n). \tag{3.20}$$

By taking the $L^2$ inner product of (3.10) with $-\frac{1}{2}U^{n+\frac{1}{2}}$, we obtain

$$-\frac{1}{4}(\|U^{n+1}\|^2 - \frac{1}{4}\|U^n\|^2) = -(\phi^{\star,n+\frac{1}{2}}(\phi^{n+1} - \phi^n), U^{n+\frac{1}{2}}) \tag{3.21}$$

By combining (3.18), (3.20) and (3.21), we obtain

$$\frac{\epsilon^2}{2}(\|\nabla \phi^{n+1}\|^2 - \|\nabla \phi^n\|^2) + \frac{1}{4}(\|U^{n+1}\|^2 - \|U^n\|^2) \\ + \alpha(\psi^{n+\frac{1}{2}}, \phi^{n+1} - \phi^n) = -\delta t M \|\nabla w^{n+\frac{1}{2}}\|^2. \tag{3.22}$$

We derive from (3.11) that

$$-\Delta(\psi^{n+1} - \psi^n) = \phi^{n+1} - \phi^n - (\overline{\phi}^{n+1} - \overline{\phi}^n). \tag{3.23}$$

By taking the $L^2$ inner product of (3.23) with $\alpha \psi^{n+\frac{1}{2}}$ and notice $\int_\Omega \psi^{n+\frac{1}{2}} d\boldsymbol{x} = 0$, we obtain

$$\frac{\alpha}{2}(\|\nabla \psi^{n+1}\|^2 - \|\nabla \psi^n\|^2) = \alpha(\psi^{n+\frac{1}{2}}, \phi^{n+1} - \phi^n) - \alpha(\overline{\phi}^{n+1} - \overline{\phi}^n)(1, \psi^{n+\frac{1}{2}}) \\ = \alpha(\psi^{n+\frac{1}{2}}, \phi^{n+1} - \phi^n). \tag{3.24}$$

By combining (3.22) and (3.24), we obtain

$$\frac{\epsilon^2}{2}(\|\nabla \phi^{n+1}\|^2 - \|\nabla \phi^n\|^2) + \frac{1}{4}(\|U^{n+1}\|^2 - \|U^n\|^2) \\ + \frac{\alpha}{2}(\|\nabla \psi^{n+1}\|^2 - \|\nabla \psi^n\|^2) = -\delta t M \|\nabla w^{n+\frac{1}{2}}\|^2. \tag{3.25}$$

Then proof is complete. $\square$

We can easily construct another second order scheme based on the backward differentiation formula.

**Scheme 2.** *Assuming that $\phi^{n-1}$, $\phi^n$ and $U^{n-1}$, $U^n$ are known, we solve $\phi^{n+1}$, $U^{n+1}$ as follows:*

$$\frac{3\phi^{n+1} - 4\phi^n + \phi^{n-1}}{2\delta t} = M\Delta w^{n+1}, \tag{3.26}$$

$$w^{n+1} = -\epsilon^2 \Delta \phi^{n+1} + \phi^{\dagger,n+1} U^{n+1} + \alpha \psi^{n+1}, \tag{3.27}$$

$$3U^{n+1} - 4U^n + U^{n-1} = 2\phi^\dagger(3\phi^{n+1} - 4\phi^n + \phi^{n-1}). \tag{3.28}$$

where $\phi^{\dagger,n+1} = 2\phi^n - \phi^{n-1}$, and

$$(3.29) \quad \begin{cases} -\Delta\psi^{n+1} = \phi^{n+1} - \overline{\phi}^{n+1}, \\ \int_\Omega \psi^{n+1} d\boldsymbol{x} = 0, \end{cases}$$

where $\overline{\phi}^{n+1} = \frac{1}{|\Omega|} \int_\Omega \phi^{n+1} d\boldsymbol{x}$.

This scheme possesses the same nice properties as the **Scheme 1**.

3.2. **PF-BCP-NS model.** We now consider the hydrodynamically coupled phase field diblock copolymer model (2.19)-(2.22). As in the above, we introduce an auxiliary function $U = \phi^2 - 1$ so that the total energy (2.18) becomes

$$(3.30) \quad E(\boldsymbol{u}, \phi, U) = \int_\Omega \Big(\frac{1}{2}|\boldsymbol{u}|^2 + \frac{\epsilon^2}{2}|\nabla\phi|^2 + \frac{1}{4}U^2 + \frac{\alpha}{2}|\nabla\psi|^2\Big) d\boldsymbol{x},$$

and the model (2.19)-(2.22) becomes

$$(3.31) \quad \phi_t = M\Delta w,$$
$$(3.32) \quad w = \lambda(-\epsilon^2 \Delta\phi + \phi U + \alpha\psi),$$
$$(3.33) \quad U_t = 2\phi\phi_t,$$
$$(3.34) \quad \boldsymbol{u}_t + (\boldsymbol{u}\cdot\nabla)\boldsymbol{u} + \nabla p - \nu\Delta\boldsymbol{u} + \phi\nabla w = 0,$$
$$(3.35) \quad \nabla\cdot\boldsymbol{u} = 0,$$

with the boundary conditions

$$(3.36) \quad (i)\ \boldsymbol{u}, p, \phi, w \text{ are periodic; or}$$
$$(3.37) \quad (ii)\ \boldsymbol{u}|_{\partial\Omega} = \partial_{\mathbf{n}}\phi|_{\partial\Omega} = \nabla w \cdot \mathbf{n}|_{\partial\Omega} = 0.$$

It is easy to see that the transformed PDE system (3.31)–(3.35) is formally equivalent to the original PF-BCP-NS model (2.19)-(2.22). By taking the $L^2$ inner product of (3.31) with $w$, of (3.32) with $-\phi_t$, of (3.33) with $\frac{1}{2}U$, of (3.34) with $\boldsymbol{u}$, performing integration by parts and summing all equalities up, we obtain the energy dissipation law of the new system (3.31)–(3.35):

$$(3.38) \quad \frac{d}{dt}E(\boldsymbol{u},\phi,U) = -M\|\nabla w\|^2 - \nu\|\nabla\boldsymbol{u}\|^2 \leq 0.$$

We now develop the second order semi-discrete numerical schemes based on the crank-Nicolson scheme to solve the flow field coupled system (3.31)-(3.35).

**Scheme 3.** *Assuming that $\boldsymbol{u}^{n-1}, \boldsymbol{u}^n, \phi^{n-1}, \phi^n, U^{n-1}, U^n, p^n$ are known, we compute the $\boldsymbol{u}^{n+1}, \phi^{n+1}, w^{n+\frac{1}{2}}, U^{n+1}, p^{n+1}$ as follows:*
**Step 1:**

$$(3.39) \quad \frac{\phi^{n+1} - \phi^n}{\delta t} + \nabla\cdot(\widetilde{\boldsymbol{u}}^{n+\frac{1}{2}} \phi^{\star,n+\frac{1}{2}}) = M\Delta w^{n+\frac{1}{2}},$$

$$(3.40) \quad w^{n+\frac{1}{2}} = \lambda(-\Delta\phi^{n+\frac{1}{2}} + \phi^{\star,n+\frac{1}{2}} U^{n+\frac{1}{2}} + \alpha\psi^{n+\frac{1}{2}}),$$

$$(3.41) \quad U^{n+1} - U^n = 2\phi^{\star,n+\frac{1}{2}}(\phi^{n+1} - \phi^n),$$

$$(3.42) \quad \frac{\widetilde{\boldsymbol{u}}^{n+1} - \boldsymbol{u}^n}{\delta t} + B(\boldsymbol{u}^{\star,n+\frac{1}{2}}, \widetilde{\boldsymbol{u}}^{n+\frac{1}{2}}) - \nu\Delta\widetilde{\boldsymbol{u}}^{n+\frac{1}{2}} + \nabla p^n + \phi^{\star,n+\frac{1}{2}}\nabla w^{n+\frac{1}{2}} = 0,$$

where $B(\boldsymbol{u}, \mathbf{v}) = (\boldsymbol{u} \cdot \nabla)\mathbf{v} + \frac{1}{2}(\nabla \cdot \boldsymbol{u})\mathbf{v}$, $\boldsymbol{u}^{\star,n+\frac{1}{2}} = \frac{3}{2}\boldsymbol{u}^n - \frac{1}{2}\boldsymbol{u}^{n-1}$, $\phi^{\star,n+\frac{1}{2}} = \frac{3}{2}\phi^n - \frac{1}{2}\phi^{n-1}$ and $\widetilde{\boldsymbol{u}}^{n+\frac{1}{2}} = \frac{\widetilde{\boldsymbol{u}}^{n+1} + \boldsymbol{u}^n}{2}$. $\psi^{n+\frac{1}{2}}$ is defined by

$$(3.43) \quad \begin{cases} -\Delta \psi^{n+\frac{1}{2}} = \phi^{n+\frac{1}{2}} - \overline{\phi}^{n+\frac{1}{2}}, \\ \int_\Omega \psi^{n+\frac{1}{2}} d\boldsymbol{x} = 0, \end{cases}$$

with $\overline{\phi}^{n+\frac{1}{2}} = \frac{1}{|\Omega|} \int_\Omega \phi^{n+\frac{1}{2}} d\boldsymbol{x}$.

**Step 2:**

$$(3.44) \quad \frac{\boldsymbol{u}^{n+1} - \widetilde{\boldsymbol{u}}^{n+1}}{\delta t} + \nabla\left(\frac{p^{n+1} - p^n}{2}\right) = 0,$$

$$(3.45) \quad \nabla \cdot \boldsymbol{u}^{n+1} = 0.$$

The boundary conditions can be either

$$(3.46) \quad (i)\ \widetilde{\boldsymbol{u}}^{n+1}, \boldsymbol{u}^{n+1}, p^{n+1}, \phi^{n+1}, w^{n+\frac{1}{2}} \ \text{are periodic; or}$$

$$(3.47) \quad (ii)\ \widetilde{\boldsymbol{u}}^{n+1}|_{\partial\Omega} = \boldsymbol{u}^{n+1} \cdot \mathbf{n}|_{\partial\Omega} = \partial_{\mathbf{n}}\phi^{n+1}|_{\partial\Omega} = \nabla w^{n+\frac{1}{2}} \cdot \mathbf{n}|_{\partial\Omega} = 0.$$

Several remarks are in order:

**Remark 3.4.**

- $B(\boldsymbol{u}, \mathbf{v})$ is the skew-symmetric form [42] of the nonlinear advection term $(\boldsymbol{u} \cdot \nabla)\boldsymbol{u}$ in the Naiver-Stokes equatiosn, If the velocity is divergence free, then $B(\boldsymbol{u}, \boldsymbol{u}) = (\boldsymbol{u} \cdot \nabla)\boldsymbol{u}$. Moreover, we have

$$(3.48) \quad (B(\boldsymbol{u}, \mathbf{v}), \mathbf{v}) = 0, \quad \forall \boldsymbol{u}, \mathbf{v} \in H^1, \text{ and } \boldsymbol{u} \cdot \mathbf{n}|_{\partial\Omega} = \mathbf{0} \text{ or } \boldsymbol{u} \text{ is periodic.}$$

- A second order pressure correction scheme [43] is used to decouple the computations of pressure from that of the velocity. This projection methods are analyzed in [32] where it is shown (discrete time, continuous space) that the schemes are second order accurate for velocity in $\ell^2(0, T; L^2(\Omega))$ but only first order accurate for pressure in $\ell^\infty(0, T; L^2(\Omega))$. The loss of accuracy for pressure is due to the artificial boundary condition (3.44) imposed on pressure [10]. We also remark that the Crank-Nicolson scheme with linear extrapolation is a popular time discretization for the Navier-Stokes equation. We refer to [17, 18, 22] and references therein for analysis on this type of discretization.

- **Step 1** consists of a coupled system for $\boldsymbol{u}^{n+1}, \phi^{n+1}, w^{n+\frac{1}{2}}$ and $U^{n+1}$. As in the last section, we can eliminate $w^{n+\frac{1}{2}}$ and $U^{n+1}$ from (3.39)-(3.41) to obtain:

$$(3.49) \quad \begin{aligned} \frac{\phi^{n+1} - \phi^n}{\delta t} &+ \nabla \cdot (\widetilde{\boldsymbol{u}}^{n+\frac{1}{2}} \phi^{\star,n+\frac{1}{2}}) \\ &= M\left(-\epsilon^2 \Delta^2 \phi^{n+\frac{1}{2}} + \Delta((\phi^{\star,n+\frac{1}{2}})^2 \phi^{n+\frac{1}{2}}) - \alpha(\phi^{n+\frac{1}{2}} - \overline{\phi}^{n+\frac{1}{2}})\right) + g^n. \end{aligned}$$

- In order to solve the coupled system (3.49)-(3.42), we first apply $(-\Delta)^{-1}$ to (3.49) to obtain

$$(3.50)\ \left(\frac{2}{M\delta t} + \alpha\right)(-\Delta)^{-1}\phi^{n+1} - \epsilon^2 \Delta \phi^{n+1} + (\phi^{\star,n+\frac{1}{2}})^2 \phi^{n+1} + \frac{1}{M}(-\Delta)^{-1}\nabla \cdot (\widetilde{\boldsymbol{u}}^{n+\frac{1}{2}} \phi^{\star,n+\frac{1}{2}}) = h^n.$$

*The coupled system (3.50)-(3.42) can then be solved by using a preconditioned conjugate gradient method with a block diagonal preconditioner associated with the following system:*

$$(\frac{2}{M\delta t} + \alpha)(-\Delta)^{-1}\phi^{n+1} - \epsilon^2 \Delta \phi^{n+1} = f_1^n, \tag{3.51}$$

$$\frac{1}{2\delta t}\widetilde{\boldsymbol{u}}^{n+1} - \nu \Delta \widetilde{\boldsymbol{u}}^{n+1} = f_2^n. \tag{3.52}$$

*It is shown in [36] that this kind of preconditioner is very effective for a simple phase-field model of two-phase incompressible flow. Our numerical experiments also show that it is effective for this PF-BCP-NS model.*

- *By taking the divergence of (3.44), we find that* **Step 2** *is equivalent to a Poisson equation for $p^{n+1}$ which can be efficiently solved by one's favorite method, in particular by spectral-Galerkin method in both periodic case or non-periodic case.*

We have the following energy stability proof.

**Theorem 3.2.** *The scheme (3.39)-(3.45) admits a unique solution. Furthermore, it is unconditionally energy stable in the sense that it satisfies the following discrete energy dissipation law:*

$$\frac{1}{\delta t}(E_{ns-cn2}(\boldsymbol{u}^{n+1}, \phi^{n+1}, U^{n+1}, p^{n+1}) - E_{ns-cn2}(\boldsymbol{u}^n, \phi^n, U^n, p^n)) = -M\|\nabla w^{n+\frac{1}{2}}\|^2 - \nu\|\nabla \widetilde{\boldsymbol{u}}^{n+\frac{1}{2}}\|^2,$$

*where $E_{ns-cn2}(\boldsymbol{u}, \phi, U, p) = \frac{1}{2}\|\boldsymbol{u}\|^2 + \lambda(\frac{\epsilon^2}{2}\|\nabla\phi\|^2 + \frac{1}{4}\|U\|^2 + \frac{\alpha}{2}\|\nabla\psi\|^2) + \frac{\delta t^2}{8}\|\nabla p\|^2$, where $\psi$ is determined by $\phi$ from (3.43).*

*Proof.* Taking the $L^2$ inner product of (3.39) with $\delta t w^{n+\frac{1}{2}}$, we obtain

$$(\phi^{n+1} - \phi^n, w^{n+\frac{1}{2}}) - \delta t(\widetilde{\boldsymbol{u}}^{n+\frac{1}{2}}\phi^{\star,n+\frac{1}{2}}, \nabla w^{n+\frac{1}{2}}) = -\delta t M\|\nabla w^{n+\frac{1}{2}}\|^2. \tag{3.53}$$

For (3.40) and (3.41), we perform the same argument as in the proof of Theorem 3.1, and combining with (3.53), we obtain

$$\begin{aligned}\lambda\frac{\epsilon^2}{2}(\|\nabla\phi^{n+1}\|^2 - \|\nabla\phi^n\|^2) + \lambda\frac{1}{4}(\|U^{n+1}\|^2 - \|U^n\|^2) - \delta t(\widetilde{\boldsymbol{u}}^{n+\frac{1}{2}}\phi^{\star,n+\frac{1}{2}}, \nabla w^{n+\frac{1}{2}}) \\ + \lambda\frac{\alpha}{2}(\|\nabla\psi^{n+1}\|^2 - \|\nabla\psi^n\|^2) = -\delta t M\|\nabla w^{n+\frac{1}{2}}\|^2.\end{aligned} \tag{3.54}$$

Taking the $L^2$ inner product of (3.42) with $\delta t \widetilde{\boldsymbol{u}}^{n+\frac{1}{2}}$, we obtain

$$\frac{1}{2}\|\widetilde{\boldsymbol{u}}^{n+1}\|^2 - \frac{1}{2}\|\boldsymbol{u}^n\|^2 + \nu\delta t\|\nabla\widetilde{\boldsymbol{u}}^{n+\frac{1}{2}}\|^2 + \delta t(\nabla p^n, \widetilde{\boldsymbol{u}}^{n+\frac{1}{2}}) + \delta t(\phi^{\star,n+\frac{1}{2}}\nabla w^{n+\frac{1}{2}}, \widetilde{\boldsymbol{u}}^{n+\frac{1}{2}}) = 0. \tag{3.55}$$

Taking the $L^2$ inner product of (3.44) with $\delta t \boldsymbol{u}^{n+1}$, we obtain

$$\frac{1}{2}\|\boldsymbol{u}^{n+1}\|^2 + \frac{1}{2}\|\boldsymbol{u}^{n+1} - \widetilde{\boldsymbol{u}}^{n+1}\|^2 = \frac{1}{2}\|\widetilde{\boldsymbol{u}}^{n+1}\|^2, \tag{3.56}$$

where we used explicitly the divergence free condition for $\boldsymbol{u}^{n+1}$.

Next, we rewrite the projection step (3.44) as

$$\boldsymbol{u}^{n+1} + \boldsymbol{u}^n - 2\widetilde{\boldsymbol{u}}^{n+\frac{1}{2}} + \frac{\delta t}{2}\nabla(p^{n+1} - p^n) = 0. \tag{3.57}$$

Taking the $L^2$ inner product of (3.57) with $\delta t \frac{1}{2}\nabla p^n$, we arrive at

$$\frac{\delta t^2}{8}(\|\nabla p^{n+1}\|^2 - \|\nabla p^n\|^2 - \|\nabla(p^{n+1} - p^n)\|^2) = \delta t(\nabla p^n, \widetilde{\boldsymbol{u}}^{n+\frac{1}{2}}). \tag{3.58}$$

Furthermore, it follows directly from (3.44) that

$$(3.59) \quad \frac{\delta t^2}{8}\|\nabla(p^{n+1} - p^n)\|^2 = \frac{1}{2}\|\boldsymbol{u}^{n+1} - \widetilde{\boldsymbol{u}}^{n+1}\|^2.$$

Combining (3.54), (3.55), (3.56), (3.58) and (3.59), we have

$$\frac{1}{2}(\|\boldsymbol{u}^{n+1}\|^2 - \|\boldsymbol{u}^n\|^2) + \frac{\delta t^2}{8}(\|\nabla p^{n+1}\|^2 - \|\nabla p^n\|^2)$$
$$+ \lambda\Big(\frac{\epsilon^2}{2}(\|\nabla \phi^{n+1}\|^2 - \|\nabla \phi^n\|^2) + \frac{1}{4}(\|U^{n+1}\|^2 - \|U^n\|^2) + \frac{\alpha}{2}(\|\nabla \psi^{n+1}\|^2 - \|\nabla \psi^n\|^2)\Big)$$
$$= -\delta t M \|\nabla w^{n+\frac{1}{2}}\|^2 - \nu \delta t \|\nabla \widetilde{\boldsymbol{u}}^{n+\frac{1}{2}}\|^2.$$

From the above, we immediately derive that if $(\boldsymbol{u}^n, \widetilde{\boldsymbol{u}}^n, \phi^n, w^n, U^n, p^n) = 0$, then we have the next step numerical sloution $(\boldsymbol{u}^{n+1}, \widetilde{\boldsymbol{u}}^{n+1}, \phi^{n+1}, w^{n+1}, U^{n+1}, p^{n+1}) = 0$. Hence, the uniqueness is proved. The existence can be established by using the above stability result coupled with a generalized Lax-Milgram theorem.

□

## 4. Numerical Simulations

We now present various numerical experiments to verify the stability and accuracy of the proposed numerical schemes. We set $\Omega = [0, 2\pi]^2$ and use the Fourier-spectral method to discretize the space variables. In all computations, we use $129 \times 129$ Fourier modes, and set

$$(4.1) \quad \epsilon = 0.06, \lambda = 1, \nu = 1.$$

### 4.1. Accuracy test.
We first perform numerical simulations to test the convergence rates of the three proposed schemes, i.e., Scheme-1 ((3.8)-(3.10)) and Scheme-2 ((3.26)-(3.28)) for the PF-BCP model, and Scheme-3 ((3.39)-(3.45)) for the PF-BCP-NS model.

#### 4.1.1. *With a manufactured exact solution.*
In the first example, we use the exact solution

$$(4.2) \quad \phi(x, y, t) = \Big(\frac{\sin(2x)\sin(2y)}{4} + 0.48\Big)\Big(1 - \frac{\sin^2(t)}{2}\Big)$$

for the PF-BCP system (2.8)-(2.9). In Table 1, we list the $L^2$ errors of the phase variable $\phi$ between the numerical solution and the exact solution at $t = 0.1$ with different time step sizes. We observe that both of the schemes Scheme-1 and Scheme-2 achieve second order accuracy in time.

For the PF-BCP-NS model (2.19)-(2.22), we take the exact solutions to be

$$(4.3) \quad \begin{cases} \phi(x, y, t) = \Big(\dfrac{\sin(2x)\sin(2y)}{4} + 0.48\Big)\Big(1 - \dfrac{\sin^2(t)}{2}\Big), \\ u(x, y, t) = \sin(2y)\sin(x)^2 \sin(t), \\ v(x, y, t) = -\sin(2x)\sin(y)^2 \sin(t), \\ p(x, y, t) = \cos(x)\sin(y)\sin(t). \end{cases}$$

In Table 2, we list the $L^2$ errors of the phase variable $\phi$, the velocity field $\boldsymbol{u} = (u, v)$ and pressure $p$ between the numerical solution and the exact solution at $t = 0.1$ with different time step sizes. We observe that Scheme-3 also achieves second order accuracy in time.

| $\delta t$ | Scheme-1 | Order | Scheme-2 | Order |
|---|---|---|---|---|
| $2 \times 10^{-2}$ | $1.43E(-4)$ | – | $1.21E(-4)$ | – |
| $1 \times 10^{-2}$ | $3.60E(-5)$ | 1.98 | $3.05E(-5)$ | 1.98 |
| $5 \times 10^{-3}$ | $9.00E(-6)$ | 2.00 | $7.68E(-6)$ | 1.98 |
| $2.5 \times 10^{-3}$ | $2.25E(-6)$ | 2.00 | $1.90E(-6)$ | 2.01 |
| $1.25 \times 10^{-3}$ | $5.63E(-7)$ | 1.99 | $4.78E(-7)$ | 1.99 |
| $6.25 \times 10^{-4}$ | $1.40E(-7)$ | 2.00 | $1.20E(-7)$ | 1.99 |
| $3.125 \times 10^{-4}$ | $3.52E(-8)$ | 1.99 | $3.04E(-8)$ | 1.98 |
| $1.5625 \times 10^{-4}$ | $8.81E(-9)$ | 1.99 | $7.63E(-9)$ | 1.99 |

TABLE 1. Accuracy test with given exact solution test for the PF-BCP model (2.4)-(2.5). The $L^2$ errors at $t = 0.1$ for the phase variable $\phi$, computed by the Scheme-1 and Scheme-2 using various time steps.

| $\delta t$ | $u$ | Order | $v$ | Order | $p$ | Order | $\phi$ | order |
|---|---|---|---|---|---|---|---|---|
| $8 \times 10^{-3}$ | $2.30E(-7)$ | – | $2.31E(-7)$ | – | $1.58E(-5)$ | – | $1.96E(-5)$ | – |
| $4 \times 10^{-3}$ | $5.20E(-8)$ | 2.14 | $5.20E(-8)$ | 2.15 | $4.09E(-6)$ | 1.94 | $4.97E(-6)$ | 1.97 |
| $2 \times 10^{-3}$ | $1.21E(-8)$ | 2.10 | $1.21E(-8)$ | 2.10 | $1.00E(-6)$ | 2.03 | $1.24E(-6)$ | 2.00 |
| $1 \times 10^{-3}$ | $2.98E(-9)$ | 2.02 | $2.98E(-9)$ | 2.02 | $2.52E(-7)$ | 1.98 | $3.12E(-7)$ | 1.99 |
| $5 \times 10^{-4}$ | $7.22E(-10)$ | 2.04 | $7.23E(-10)$ | 2.04 | $6.34E(-8)$ | 1.99 | $7.79E(-8)$ | 2.00 |

TABLE 2. Accuracy test with given exact solution for PF-BCP-NS model (2.19)-(2.22). The $L^2$ errors at $t = 0.1$ for the phase variable $\boldsymbol{u} = (u, v)$, $p$ and $\phi$, computed by Scheme-3 using various time steps.

4.1.2. *With a given initial condition.* Next, we further examine the temporal accuracy using an example with a given initial condition but with no explicit exact solution. The initial condition of $\phi$ is set to be

(4.4) $$\phi(x, y, 0) = (\frac{\sin(2x)\sin(2y)}{4} + 0.48).$$

Since the exact solutions are not known, we choose the solution obtained by Scheme-1 with the time step size $\delta t = 1 \times 10^{-5}$ as the benchmark solution for computing errors. We present the $L^2$ errors of the phase variable between the numerical solution and the benchmark solution at $t = 1$ with different time step sizes in Table 3. We observe that both Scheme-1 and Scheme-2 are second order accurate in time for this example.

| $\delta t$ | Scheme-1 | Order | Scheme-2 | Order |
|---:|:---:|:---:|:---:|:---:|
| $2 \times 10^{-2}$ | $3.33E(-4)$ | – | $3.35E(-4)$ | – |
| $1 \times 10^{-2}$ | $1.01E(-4)$ | 1.72 | $1.01E(-4)$ | 1.73 |
| $5 \times 10^{-3}$ | $2.83E(-5)$ | 1.83 | $2.83E(-5)$ | 1.83 |
| $2.5 \times 10^{-3}$ | $7.52E(-6)$ | 1.91 | $7.53E(-6)$ | 1.91 |
| $1.25 \times 10^{-3}$ | $1.94E(-6)$ | 1.95 | $1.94E(-6)$ | 1.96 |
| $6.25 \times 10^{-4}$ | $4.93E(-7)$ | 1.97 | $4.94E(-7)$ | 1.97 |
| $3.125 \times 10^{-4}$ | $1.24E(-7)$ | 1.99 | $1.27E(-7)$ | 1.95 |
| $1.5625 \times 10^{-4}$ | $3.11E(-8)$ | 1.99 | $3.11E(-8)$ | 2.02 |

TABLE 3. Accuracy test of Scheme-1 and Scheme-2 with given initial condition for the PF-BCP model (2.4)-(2.5). The $L^2$ errors at $t = 1$ for the phase variable $\phi$.

4.2. **Phase separation.** In this example, we study the molecular self-assembly in BCP thin films to form lamellar or cylindrical nanostructures through phase separation, or referred as to spinodal decomposition. The process of the phase separation can be studied by considering a homogeneous binary mixture, which is quenched into the unstable part of its miscibility gap. In this case, the spinodal decomposition takes place, which manifests in the spontaneous growth of the concentration fluctuations that leads the system from the homogeneous state to the two-phase state. Shortly after the phase separation starts, the domains of the binary components are formed and the interface between the two phases can be specified [2, 9, 57].

The initial conditions are taken as a randomly perturbed concentration field as follows:
$$\phi(t=0) = \widehat{\phi}_0 + 0.01\,\mathrm{rand}(x,y), \tag{4.5}$$
where the $\mathrm{rand}(x,y)$ is a uniformly distributed random function in $[-1,1]^2$ with zero mean. We present below numerical simulations for the PF-BCP model (2.4)-(2.5), PF-BCP-NS model (2.19)-(2.22), and PF-BCP model with imposed electric field (4.6)-(4.7) by varying the nonlocal parameter $\alpha$ and initial average value $\widehat{\phi}_0$.

4.2.1. *For PF-BCP model.* We first choose a very small nonlocal parameter $\alpha = 0.001$. For this choice, the PF-BCP model (2.4)-(2.5) is expected to be consistent with the standard Cahn-Hilliard model [5, 14, 25, 27].

Before we run the phase separation simulation to the steady state, we need to choose a suitable time step. Even though any time step size $\delta t$ is allowable for the computations from the stability concern since all developed schemes are unconditionally energy stable, we emphasize that larger time step will definitely induce large numerical errors. Therefore, we need to discover the rough range of the allowable maximum time step size in order to obtain good accuracy and to consume as low computational cost as possible. This time step range could be estimated through the energy evolution curve plots, shown in Fig. 1, where we compare the time evolution of the free energy for five different time step sizes until $t = 1$ using Scheme-1. We observe that all five energy curves show decays monotonically for all time step sizes, which numerically confirms that our algorithms are unconditionally energy stable. For smaller time steps of $\delta t = 0.0001, 0.0005, 0.001, 0.005$, the four energy curves coincide very well. But for the larger time step of $\delta t = 0.01$ ($\sim O(\epsilon)$), the energy curve deviates viewable away from others. This means the adopted time step size should not be larger than 0.01, in order to get reasonably good accuracy. Thus we choose $\delta t = 0.001$ for all simulations.

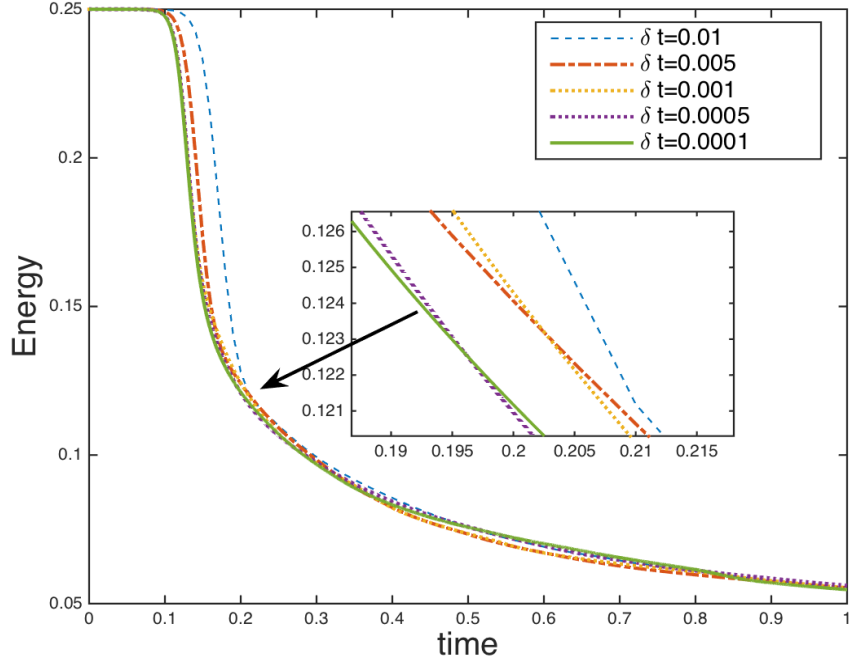

Figure 1. Time evolution of the free energy functional for five different time steps of $\delta t = 0.0001, 0.0005, 0.001, 0.005,$ and $0.01$ for $\widehat{\phi}_0 = 0$. The energy curves show the decays for all time steps, which confirms that our algorithm is unconditionally energy stable. The small inset figure shows the small differences in the energy evolution for all four time steps.

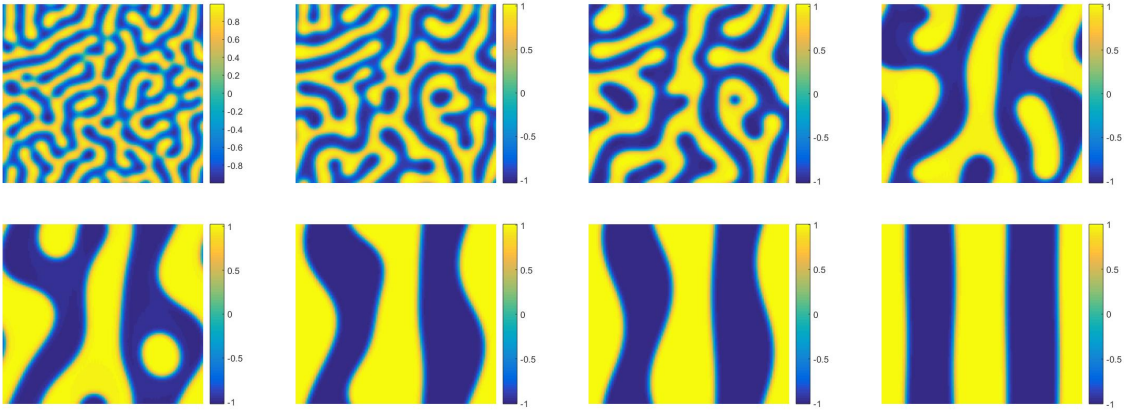

Figure 2. The 2D dynamical evolution of the phase variable $\phi$ for the PF-BCP model with the initial condition $\widehat{\phi}_0 = 0, \alpha = 0.001$ and time step $\delta t = 0.001$. Snapshots of the numerical approximation are taken at $t = 0.25, 0.5, 1, 5, 10, 20, 30, 100$.

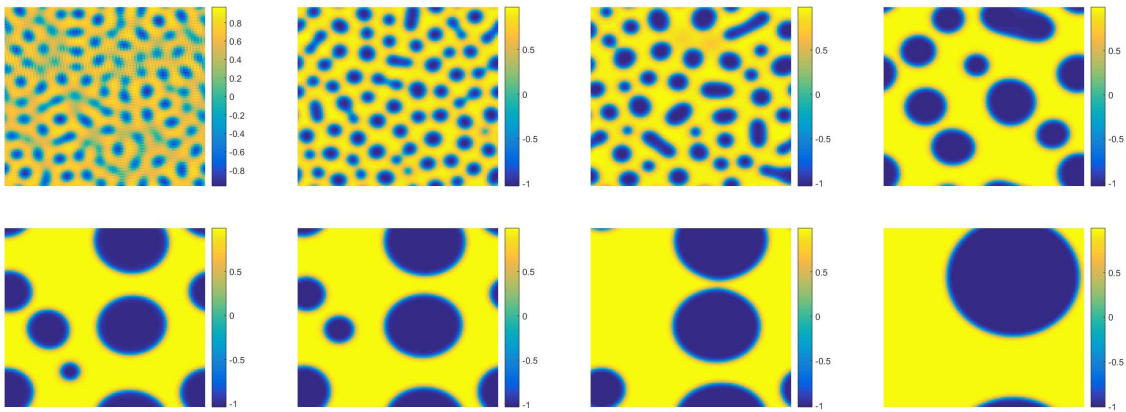

FIGURE 3. The 2D dynamical evolution of the phase variable $\phi$ for the PF-BCP model with the initial condition $\widehat{\phi}_0 = 0.3, \alpha = 0.001$ and time step $\delta t = 0.001$. Snapshots of the numerical approximation are taken at $t = 0.25, 0.5, 1, 10, 40, 60, 100, 400$.

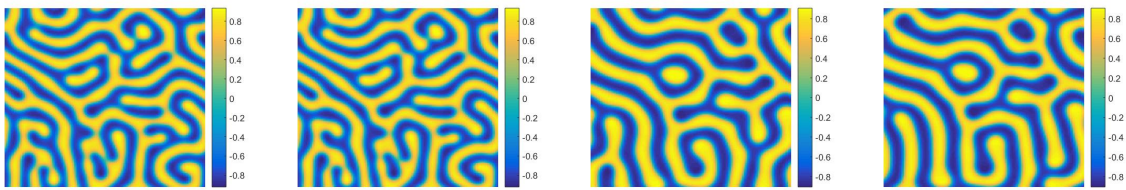

FIGURE 4. The 2D dynamical evolution of the phase variable $\phi$ for the PF-BCP model with the initial condition $\widehat{\phi}_0 = 0, \alpha = 5$ and time step $\delta t = 0.001$. Snapshots of the numerical approximation are taken at $t = 0.25, 0.5, 40, 700$.

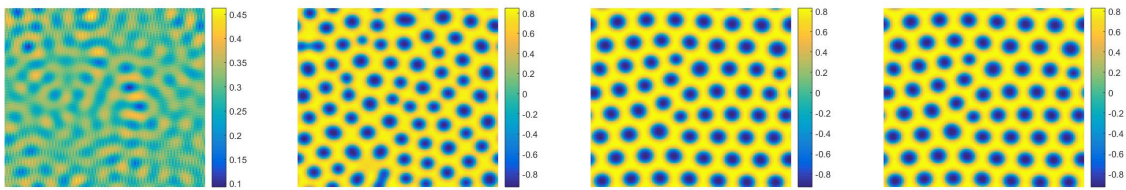

FIGURE 5. The 2D dynamical evolution of the phase variable $\phi$ for the PF-BCP model with the initial condition $\widehat{\phi}_0 = 0.3, \alpha = 10$ and time step $\delta t = 0.001$. Snapshots of the numerical approximation are taken at $t = 0.25, 1, 60, 700$.

In Fig. 2, we perform numerical simulations for initial values $\widehat{\phi}_0 = 0$, that means the volume of two monomer blocks is almost identical. We observe that phase dislocations immediately appear in the very beginning during the process of phase coarsening ($t = 0.25$ for instance). The final equilibrium solution is obtained after $t = 100$, where the lamellar (banded) nanostructures finally formed. In Fig. 3, we take the initial value of $\widehat{\phi}_0 = 0.3$ where the volume of one monomer (yellow region of $\phi = 1$) is more than the volume of the other monomer (blue region of $\phi = -1$). The dynamics are quite different from the case of $\widehat{\phi}_0 = 0$. The monomers with less volume merge and finally accumulate to circular pattern. The banded shape, as well as the circular pattern are consistent with the experimental results of lamellar phase and the cylindrical state (cf. [15, 48, 53]).

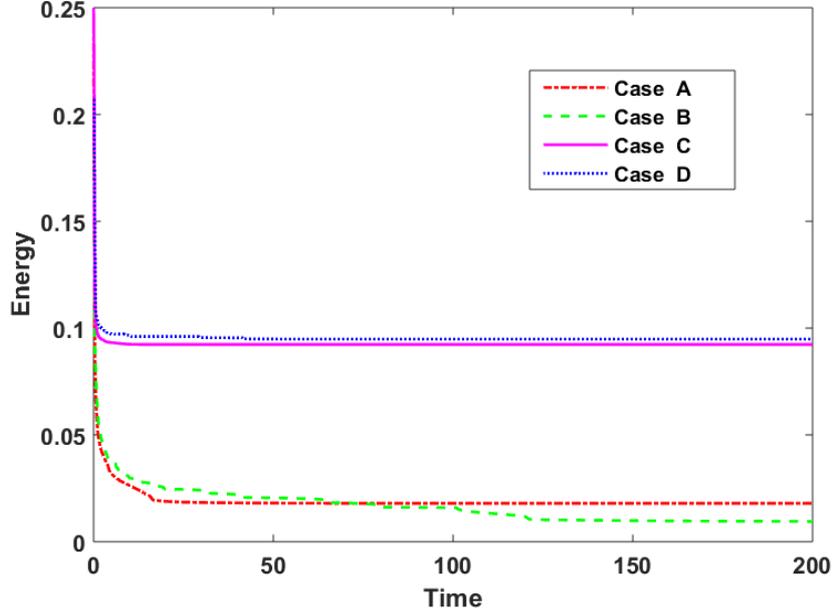

FIGURE 6. Time evolution of the free energy functional for PF-BCP model. The parameters are $(\alpha, \widehat{\phi}_0) = (0.001, 0)$ (Case A), $(\alpha, \widehat{\phi}_0) = (0.001, 0.3)$ (Case B), $(\alpha, \widehat{\phi}_0) = (5, 0)$ (Case C), $(\alpha, \widehat{\phi}_0) = (5, 0.3)$ (Case D). The energy curves decay for all time steps, which confirms that our algorithm is unconditionally energy stable.

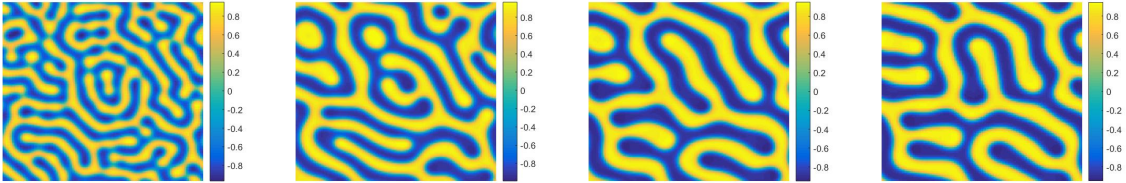

FIGURE 7. The 2D dynamical evolution of the phase variable $\phi$ for the PF-BCP-NS model with nonlocal parameter $\alpha = 5$ and the initial condition $\widehat{\phi}_0 = 0$ and time step $\delta t = 0.001$. Snapshots of the numerical approximation are taken at $t = 0.25, 1, 20, 200$.

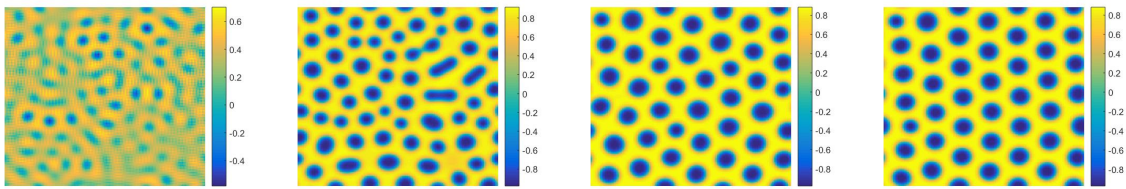

FIGURE 8. The 2D dynamical evolution of the phase variable $\phi$ for the PF-BCP-NS model with nonlocal parameter $\alpha = 5$ and the initial condition $\widehat{\phi}_0 = 0.3$ and time step $\delta t = 0.001$. Snapshots of the numerical approximation are taken at $t = 0.25, 1, 20, 200$.

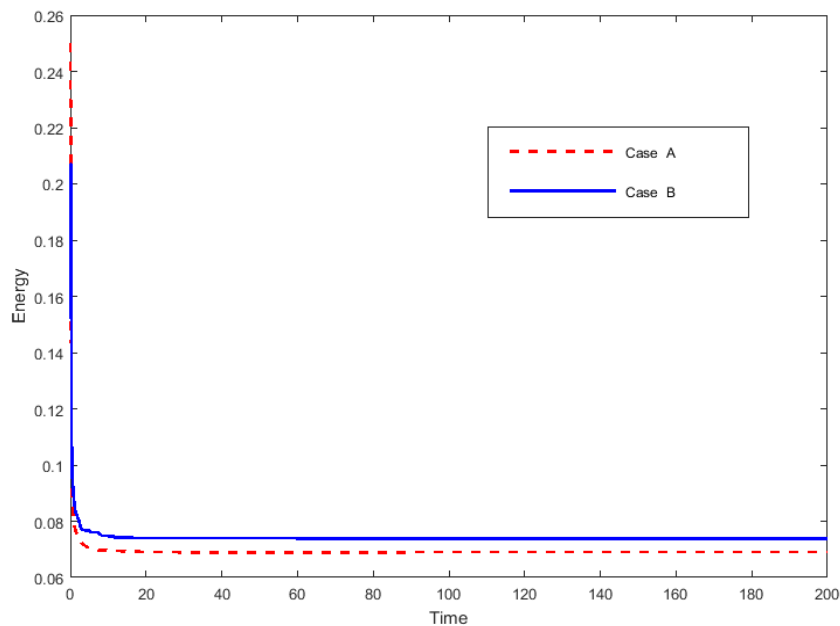

FIGURE 9. Time evolution of the free energy functional for PF-BCP-NS model for $\alpha = 5$ with various initial values of $\widehat{\phi}_0 = 0$ (Case A) and $\widehat{\phi}_0 = 0.3$ (Case B). The energy curves show the decays for all time steps, which confirms that our algorithm is unconditionally stable.

Next, we increase the nonlocal parameter to $\alpha = 5$. In Fig. 4, with the initial value of $\widehat{\phi}_0 = 0$, we observe that the initial phase separation behaviors are similar to the case of $\alpha = 0.001$ in Fig. 2. But the final lamellar phase presents very different pattern and the equilibrium is obtained at $t = 40$, that is much faster than the case of $\alpha = 0.001$.

We present in Fig. 5 results with the initial value $\widehat{\phi}_0 = 0.3$. The phase separation dynamics is still consistent to the case of $\alpha = 0.001$ in Fig. 3. The final steady state exhibits a circular pattern (cylindrical state) but with much smaller radius.

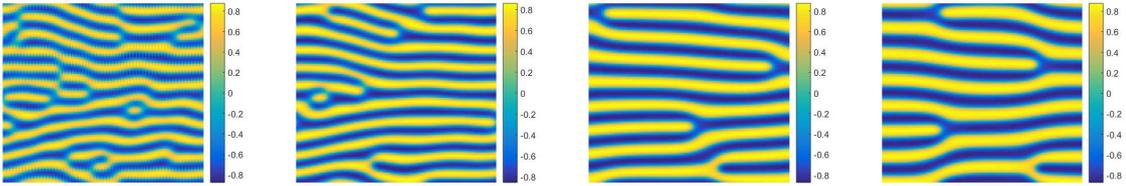

FIGURE 10. The 2D dynamical evolution of the phase variable $\phi$ for the PF-BCP model with imposed electric filed for the initial condition $\widehat{\phi}_0 = 0$ and the nonlocal parameter $\alpha = 10$ and time step $\delta t = 0.001$. Snapshots of the numerical approximation are taken at $t = 0.25, 0.5, 5, 700$.

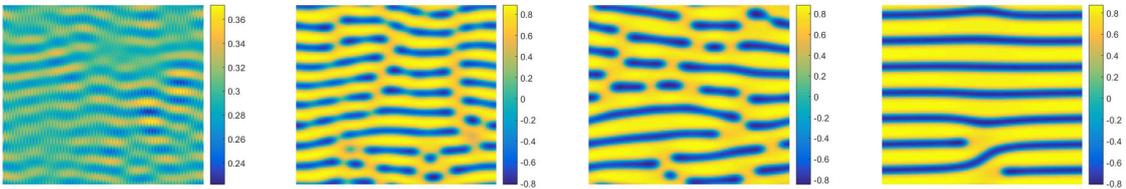

FIGURE 11. The 2D dynamical evolution of the phase variable $\phi$ for the PF-BCP model with imposed electric filed for the initial condition $\widehat{\phi}_0 = 0.3$ and the nonlocal parameter $\alpha = 10$ and time step $\delta t = 0.001$. Snapshots of the numerical approximation are taken at $t = 0.25, 0.5, 4, 700$.

In Fig. 6, we present the evolution of the free energy functional for both initial averages of $\widehat{\phi}_0 = 0, 0.3$ and $\alpha = 0, 5$. The energy curves show the decays with time that confirms that our algorithms are unconditionally energy stable.

4.2.2. *For PF-BCP-NS model.* We study next the PF-BCP-NS model (2.19)-(2.22). The nonlocal parameter is set to be $\alpha = 5$. First, we take the initial value $\widehat{\phi}_0 = 0$. In Fig. 7, we observe that the initial phase separation process is similar to the no flow case of Fig. 4. The result with the initial value of $\widehat{\phi}_0 = 0.3$ is presented in Fig. 8. The morphology behaviors and final equilibrium solutions are quite similar to the no flow case of Fig. 5. In Fig. 9, we present the evolution of the free energy functional of $\widehat{\phi}_0 = 0$ and $0.3$. The energy curves show the decays with time that confirms that our algorithms are unconditionally stable.

4.2.3. *For PF-BCP model with imposed electric field.* Applying an external electric field is one of the most efficient approach to control and produce various nano-structured materials, and has attracted substantial attentions since the original work in [46]. When an external electric filed is applied along some direction, for instance, $x-$ axis, an additional term contributed by the electric field is added to the model system (2.8)-(2.9). The new phase equation model reads as follows (cf. [44, 46]):

(4.6) $$\phi_t = M\Delta w + \beta \phi_{xx},$$
(4.7) $$w = -\epsilon^2 \Delta \phi + f(\phi) + \alpha \psi,$$

where $\beta$ is the magnitude of the electric field.

The electric field term can be viewed as an imposed external force. Note that the system (4.6)-(4.7) does not follow the energy dissipation law when $\beta \neq 0$. To solve the equation (4.6), notice that the imposed term is linear, therefore we can simply modify (3.8) in Scheme-1 as

$$(4.8) \qquad \frac{\phi^{n+1} - \phi^n}{\delta t} = M\Delta w^{n+\frac{1}{2}} + \beta \phi_{xx}^{n+\frac{1}{2}},$$

and the rest of the scheme ((3.9)-(3.10)) is still the same.

We let $\beta = 0.2$ for the following two simulations. First, we set $\widehat{\phi}_0 = 0$, in Fig. 10, we observe that phase dislocations exist in the lamellar nano-structures during the initial process of phase coarsening ($t = 0.25, 0.5$). The dislocations gradually disappear through climbing and molecular diffusion under the action of the electric driving force. The final equilibrium solution that presents lamellar phase with a few dislocations is obtained around $t = 20$. Note the convergence speed to the steady state is much faster than the case without the electric field ($t = 200$ in Fig. 2).

We further set $\widehat{\phi}_0 = 0.3$, since the volume of one monomer is less than the volume of the other one, we observe more phase dislocations during the initial process of phase coarsening ($t = 0.5, 4$ for instance) in Fig. 11. The final equilibrium solution also exhibits a lamellar phase due to the electric effects. Note for the case without the electric field, the steady state when $\widehat{\phi}_0 = 0.3$ is always the cylindrical state (cf. Fig. 3, Fig. 5 or Fig. 8). This means the electric field plays a dominant role in forming lamellar nano-structures. The obtained simulations are qualitatively consistent with the numerical results in [23, 44, 46].

## 5. Concluding remarks

In this paper, we presented a set of efficient time discretization schemes for solving the PF-BCP model and PF-BCP-NS model. The schemes are (i) second order accurate in time; (ii) unconditional energy stable; and (iii) linear and easy to implement as one only needs to solve a linear system with symmetric positive definite operator at each time step. Various numerical results are presented to validate the accuracy of our schemes. We have also presented numerical simulations to show the morphological evolutions of PF-BCP model and PF-BCP-NS model. In particular, phase separations for different cases with or without flow and with or without external electric field are investigated.